\documentclass[12pt,leqno]{article}
\usepackage{amssymb}
\usepackage[mathscr]{eucal}
\usepackage{amsmath,amssymb,latexsym,theorem,bbm}
\usepackage{color}

\newcommand{\SC}{\scriptstyle}

\newcommand{\CC}{\mathsf{C}}
\newcommand{\DD}{\mathsf{D}}

\newcommand{\NN}{\mathbb{N}}

\newcommand{\RR}{\mathbb{R}}

\newcommand{\ZZ}{\mathbb{Z}}

\newcommand{\bx}{{\boldsymbol{x}}}

\newcommand{\cA}{{\mathcal A}}
\newcommand{\cB}{{\mathcal B}}

\newcommand{\cD}{{\mathcal D}}

\newcommand{\cF}{{\mathcal F}}

\newcommand{\cM}{{\mathcal M}}

\newcommand{\cN}{{\mathcal N}}

\newcommand{\cP}{{\mathcal P}}

\newcommand{\cU}{{\mathcal U}}

\newcommand{\cX}{{\mathcal X}}

\newcommand{\cY}{{\mathcal Y}}

\newcommand{\cW}{{\mathcal W}}

\newcommand{\dd}{\mathrm{d}}

\newcommand{\slu}{{\SC\mathrm{lu}}}

\newcommand{\INARp}{\textup{INAR($p$)}}

\newcommand{\EE}{\operatorname{\mathbb{E}}}
\newcommand{\PP}{\operatorname{\mathbb{P}}}
\newcommand{\OO}{\operatorname{O}}

\newcommand{\var}{\operatorname{Var}}

\newcommand{\tS}{\widetilde{S}}

\newcommand{\vare}{\varepsilon}

\renewcommand{\mid}{\,|\,}
\newcommand{\bmid}{\,\big|\,}

\renewcommand{\leq}{\leqslant}

\newcommand{\stoch}{\stackrel{\PP}{\longrightarrow}}
\newcommand{\distr}{\stackrel{\cD}{\longrightarrow}}

\newcommand{\distre}{\stackrel{\cD}{=}}

\newcommand{\lu}{\stackrel{\slu}{\longrightarrow}}
\newcommand{\as}{\stackrel{{\mathrm{a.s.}}}{\longrightarrow}}

\newcommand{\bbone}{\mathbbm{1}}
\newcommand{\ns}{{\lfloor ns\rfloor}}
\newcommand{\nt}{{\lfloor nt\rfloor}}
\newcommand{\nT}{{\lfloor nT\rfloor}}
\newcommand{\proofend}{\hfill\mbox{$\Box$}}

\numberwithin{equation}{section}

\theoremstyle{change} \theorembodyfont{\em}
\newtheorem{Lem}{Lemma.}[section]
\newtheorem{Thm}{Theorem.}[section]

\theorembodyfont{\rm}
\newtheorem{Rem}{Remark.}[section]

\begin{document}

\begin{center}
 {\bfseries\Large A note on asymptotic behavior of critical}\\[1mm]
 {\bfseries\Large Galton-Watson processes with immigration}

\vspace*{3mm}

 {\sc\large
  M\'aty\'as $\text{Barczy}^{*,\diamond}$,
  \ D\'aniel $\text{Bezd\'any}^{**}$,
  \ \framebox[1.1\width]{Gyula $\text{Pap}$}}

\end{center}

\vskip0.2cm

\noindent
 * MTA-SZTE Analysis and Stochastics Research Group,
   Bolyai Institute, University of Szeged,
   Aradi v\'ertan\'uk tere 1, H--6720 Szeged, Hungary.

\noindent
 ** Former master student of Bolyai Institute, University of Szeged,
    Aradi v\'ertan\'uk tere 1, H-6720 Szeged, Hungary.

\noindent e-mails: barczy@math.u-szeged.hu (M. Barczy),
                  bezd.dani@gmail.com (D. Bezd\'any).

\noindent $\diamond$ Corresponding author.

\renewcommand{\thefootnote}{}
\footnote{This research is supported by the grant NKFIH-1279-2/2020 of the Ministry for Innovation and Technology, Hungary.}


\begin{abstract}
In this somewhat didactic note we give a detailed alternative proof of the known result due to Wei and Winnicki (1989) which states that under second order moment assumptions
 on the offspring and immigration distributions the sequence of appropriately scaled random step functions formed from a critical Galton-Watson process
 with immigration (starting from not necessarily zero) converges weakly towards a squared Bessel process.
The proof of Wei and Winnicki (1989) is based on infinitesimal generators,
 while we use limit theorems for random step processes towards a diffusion process due to Isp\'any and Pap (2010).
This technique was already used in Isp\'any (2008), where he proved
 functional limit theorems for a sequence of some appropriately normalized nearly critical Galton-Watson
 processes with immigration starting from zero, where the offspring means tend to its critical value 1.
As a special case of Theorem 2.1 in Isp\'any (2008) one can get back the result of Wei and Winnicki (1989) in the case of zero initial value.
In the present note we handle non-zero initial values with the technique used in Isp\'any (2008), and further, we simplify some of the arguments
 in the proof of Theorem 2.1 in Isp\'any (2008) as well.
\end{abstract}

\section{Introduction and results}
\label{section_intro_results}

The study of the limit behaviour of Galton-Watson processes has a long tradition and history, see, e.g., the famous book of Athreya and Ney \cite{AN}.
A Galton-Watson process with or without immigration is called subcritical, critical and supercritical if the mean of its offspring distribution
 is less than 1, equal to 1 and greater than 1, respectively (for more details, see later on).
For a sequence of critical Galton-Watson processes without immigration, with the same offspring distribution having finite second moment and
 with initial value independent of the offspring variables such that the initial value of the \ $n^{\mathrm{th}}$ \ branching process
 in question divided by \ $n$ \ converges in distribution as \ $n\to\infty$, \ Feller \cite{Fel} proved that the sequence of appropriately scaled random step functions converges in distribution to a non-negative diffusion process without drift (for a detailed proof based on infinitesimal generators,
 see also Ethier and Kurtz \cite[Theorem 9.1.3]{EthKur}).
Grimvall \cite[Theorem 4.4]{Gri} proved a fluctuation-type limit theorem for a sequence of nearly critical Galton-Watson processes
 without immigration:
 shifting each branching process in question by its own (deterministic) initial value, under some Lindeberg-type condition on the offspring distribution
 it was shown that the sequence of appropriately scaled random step functions formed from the sequence of shifted branching processes
 converges weakly to a Wiener process with some drift and variance
 depending on the limiting behaviour of the offspring mean and variance, respectively.
In fact, Grimvall \cite[Theorem 4.4]{Gri} generalized the corresponding result of Lindvall \cite[Theorem 1]{Lin} for a sequence of critical
  Galton-Watson processes without immigration.

In this somewhat didactic note we will focus on asymptotic behaviour of critical Galton-Watson processes with immigration.
We give a detailed alternative proof of the known result due to Wei and Winnicki \cite[Theorem 2.1]{WW} which
 states that under second order moment assumptions on the offspring and immigration distributions
 the sequence of appropriately scaled random step functions formed from a critical Galton-Watson process
 with immigration (starting from not necessarily zero) converges weakly towards a squared Bessel process, see Theorem \ref{Thm_critical}.
For historical fidelity, we mention that the convergence of finite-dimensional distributions of a sequence of Galton-Watson processes with immigration towards a continuous state and
 continuous time branching process was already studied by Kawazu and Watanabe \cite{KawWat} and Aliev \cite{Ali}.
Wei and Winnicki \cite{WW} used infinitesimal generators in their proof by referring to several results of Ethier and Kurtz \cite{EthKur},
 while we will use limit theorems for random step processes towards a diffusion process due to Isp\'any and Pap \cite{IspPap}.
This technique was already used in Isp\'any \cite{Isp}, where he proved functional limit theorems for a sequence of some appropriately normalized
 nearly critical Galton-Watson processes with immigration starting from zero, where the offspring means tend to its critical value 1 under
 some conditions on the variances of the offspring and immigration distributions.
In the present note we will handle non-zero initial values with the technique used in Isp\'any (2008), and further, we can also simplify some of the arguments
 in the proof of Theorem 2.1 in Isp\'any \cite{Isp}
 mainly due to the fact that we consider only a single critical Galton-Watson process with immigration instead of a sequence of nearly critical ones.
In Remark \ref{Rem_Isp_comparison} one can find a detailed comparison of our proof of Theorem \ref{Thm_critical} and the proof of Theorem 2.1 in Isp\'any \cite{Isp}.
 Remark \ref{Rem_initial_value} is devoted to a discussion on the role of the initial value.

We also remark that, using the technique of infinitesimal generators, Sriram \cite[Theorem 3.1]{Sri}, Isp\'any et al. \cite[Theorem 2.1]{IspPapZui} and
 Khusanbaev \cite{Khu} proved functional limit theorems for a sequence of some appropriately normalized nearly critical Galton-Watson processes
 with immigration starting from zero.
Lebedev \cite{Leb} proved the result of Sriram \cite[Theorem 3.1]{Sri} independently as well.
Li \cite{Li} provided a set of sufficient conditions for the weak convergence of a sequence of Galton-Watson processes with immigration
 to a given continuous state and continuous time branching process with immigration.
Using martingale limit theorems based on Jacod and Shiryaev \cite{JacShi},
Rahimov \cite{Rah} proved functional limit theorems for a sequence of critical Galton-Watson processes
 with generation dependent immigrations starting from zero such that the means of immigration distributions tend to infinity as the number of generation goes to infinity.

Let \ $\ZZ_+$, \ $\NN$, \ $\RR$, \ $\RR_+$  \ and \ $\RR_{++}$ \ denote the set
 of non-negative integers, positive integers, real numbers, non-negative real
 numbers and positive real numbers, respectively.
For sequences \ $(a_k)_{k\in\NN}$ \ and \ $(b_k)_{k\in\NN}$, \ where \ $b_k\in\RR_{++}$, \ $k\in\NN$, \
 the notation \ $a_k = \OO(b_k)$, \ $k\in\NN$, \ means that there
 exists a constant \ $C\in\RR_{++}$ \ such that \ $\vert a_k\vert\leq C b_k$, \ $k\in\NN$.
\ In the proofs we frequently use that for any \ $\gamma \in\RR_{++}$, \ we have \ $\sum_{k=0}^n k^\gamma = \OO(n^{\gamma+1})$ \ for \ $n\in\NN$, \
 following from \ $\sum_{k=0}^n k^\gamma \leq \int_0^{n+1} x^\gamma \,\dd x = \frac{(n+1)^{\gamma+1}}{\gamma+1} \leq \frac{2^{\gamma+1}}{\gamma+1} n^{\gamma+1}$,
 \ $n\in\NN$.
\ For a function \ $f:\RR\to\RR$, \ its positive part will be denoted by \ $f^+$.
\ Every random variable will be defined on a fixed probability space
 \ $(\Omega, \cA, \PP)$.
\ Convergence in probability is denoted by \ $\stoch$.
\ For other notations, such as equality in distribution \ $\distre$ \  and convergence in distribution \ $\distr$, \
 see the beginning of Appendix \ref{CMT}.

First we recall (single-type) Galton-Watson processes with immigration.
For each \ $k \in \ZZ_+$, \ the number of individuals in the \ $k^\mathrm{th}$ \ generation is denoted by \ $X_k$.
\ By \ $\xi_{k,j}$ \ we denote the number of the offsprings produced by the \ $j^\mathrm{th}$ \ individual
  belonging to the \ $(k-1)^\mathrm{th}$ \ generation.
The number of immigrants in the \ $k^\mathrm{th}$ \ generation will be denoted by \ $\vare_k$.
Then we have
 \begin{equation}\label{MBPI(d)}
  X_k = \sum_{j=1}^{X_{k-1}} \xi_{k,j} + \vare_k , \qquad
  k \in \NN ,
 \end{equation}
 where we define \ $\sum_{j=1}^0:=0$.
\ Here
 \ $\big\{X_0, \xi_{k,j}, \, \vare_k  : k, j \in \NN \big\}$
 \ are supposed to be independent \ $\ZZ_+$-valued random variables.
Moreover, \ $\big\{\xi_{k,j} : k, j \in \NN\big\}$ \ and \ $\{\vare_k : k \in \NN\}$ \ are
 supposed to consist of identically distributed random variables, respectively.
For notational convenience, let \ $\xi$ \ and \ $\vare$ \ be random variables such that \ $\xi \distre \xi_{1,1}$ \
 and \ $\vare \distre \vare_1$.

We suppose that \ $\EE(X_0^2)<\infty$, \ $\EE(\xi^2) < \infty$ \ and \ $\EE (\vare^2) < \infty$.
\ Introduce the notations
 \begin{gather*}
  m_\xi := \EE(\xi),\qquad
  m_\vare:= \EE(\vare), \qquad
  \sigma_\xi^2 := \var(\xi), \qquad
  \sigma_\vare^2 := \var(\vare).
 \end{gather*}
For \ $k \in \ZZ_+$, \ let \ $\cF_k^X := \sigma(X_0,X_1 , \dots, X_k)$.
\ By \eqref{MBPI(d)}, \ $\EE(X_k \mid \cF_{k-1}^X) = m_\xi X_{k-1}  + m_\vare$, \ $k\in\NN$.
\ Consequently, \ $\EE(X_k) = m_\xi \EE(X_{k-1}) + m_\vare$, \ $k \in \NN$, \ which implies
 \begin{equation}\label{EXk}
  \EE(X_k) = \EE(X_0) m_\xi^k + m_\vare \sum_{j=0}^{k-1} m_\xi^j
           = \begin{cases}
              \EE(X_0) m_\xi^k + m_\vare \frac{m_\xi^k-1}{m_\xi-1} & \text{if \ $m_\xi\ne 1$,}\\[1mm]
              \EE(X_0) + m_\vare k & \text{if \ $m_\xi = 1$,}
             \end{cases} , \qquad k \in \NN .
 \end{equation}
Hence the offspring mean \ $m_\xi$ \ plays a crucial role in the asymptotic behavior of the sequence \ $(\EE(X_k) )_{k\in\ZZ_+}$.
\ A  Galton-Watson process \ $(X_k)_{k\in\ZZ_+}$ \ with immigration is referred to respectively
 as \emph{subcritical}, \emph{critical} or \emph{supercritical} if \ $m_\xi < 1$, \ $m_\xi = 1$
 \ or \ $m_\xi > 1$ \ (see, e.g., Athreya and Ney \cite[V.3]{AN}).

We give a detailed alternative proof of the following known result due to Wei and Winnicki \cite[Theorem 2.1]{WW}
 (under an additional second order moment condition on the initial value \ $X_0$, \ which is not supposed in \cite{WW}).

\begin{Thm}[Wei and Winnicki \cite{WW}] \label{Thm_critical}
Let \ $(X_k)_{k\in\ZZ_+}$ \ be a critical Galton-Watson process with immigration such that
 \ $\EE(X_0^2) < \infty$, \ $\EE(\xi^2) < \infty$ \ and \ $\EE(\vare^2) < \infty$.
\ Then
 \begin{equation}\label{convX}
  (n^{-1} X_\nt)_{t\in\RR_+} \distr (\cX_t)_{t\in\RR_+} \qquad
  \text{as \ $n \to \infty$,}
 \end{equation}
 where the limit process \ $(\cX_t)_{t\in\RR_+}$ \ is the pathwise unique strong solution of the stochastic
 differential equation (SDE)
 \begin{equation}\label{SDE_single}
  \dd \cX_t
  = m_\vare \, \dd t
    + \sqrt{ \sigma_\xi^2 \, \cX_t^+} \, \dd \cW_t , \qquad t\in\RR_+,
 \end{equation}
 with initial value \ $\cX_0 = 0$, \ where \ $(\cW_t)_{t\in\RR_+}$ \ is a standard Wiener process.
\end{Thm}

The SDE \eqref{SDE_single} has a pathwise unique strong solution \ $(\cX_t^{(x)})_{t\in\RR_+}$ \ for all initial values \ $\cX_0^{(x)} = x \in \RR$,
 \ and if \ $x \in \RR_+$, \ then \ $\cX_t^{(x)} \in \RR_+$ \ almost surely for all \ $t \in \RR_+$, \ since \ $m_\vare, \sigma_\xi^2 \in \RR_+$, \ see,
 e.g., Ikeda and Watanabe \cite[Chapter IV, Example 8.2]{IkeWat}.
The process \ $(\cX_t^{(x)})_{t\in\RR_+}$ \ is called a squared Bessel process.

\begin{Rem}
\noindent (i) Under the conditions of Theorem \ref{Thm_critical}, we have
 \[
     \Big(  n^{-1}(X_{\nt} - \EE(X_{\nt})) \Big)_{t\in\RR_+} \distr (\cM_t)_{t\in\RR_+}\qquad
                                                          \text{as \ $n \to \infty$,}
 \]
 where
 the limit process \ $(\cM_t)_{t\in\RR_+}$ \ is the pathwise unique strong solution of the SDE
 \begin{equation*}
  \dd \cM_t
  = \sqrt{\sigma_\xi^2 (\cM_t + m_\vare t)^+ } \, \dd \cW_t , \qquad
  \cM_0 = 0 ,
 \end{equation*}
 where \ $(\cW_t)_{t\in\RR_+}$ \ is a standard Wiener process.
Indeed, by the proof of Theorem \ref{Thm_critical} (see \eqref{conv_M}), we have
 \[
     \Big(  n^{-1}(X_{\nt} - \nt m_\vare) \Big)_{t\in\RR_+} \distr (\cM_t)_{t\in\RR_+}\qquad
                                                          \text{as \ $n \to \infty$,}
 \]
 and, by \eqref{mean_X},
 \[
   n^{-1}( X_{\nt} - \EE(X_{\nt}) ) = n^{-1}( X_{\nt} - \nt m_\vare)  - n^{-1}\EE(X_0), \qquad n\in\NN,\quad t\in\RR_+.
 \]

\noindent (ii) Under the conditions of Theorem \ref{Thm_critical}, in the special case of \ $\sigma_\xi=0$, \ we have \ $\PP(\xi=1)=1$, \
 and \ $(n^{-1} X_\nt)_{t\in\RR_+}\distr (m_\vare t)_{t\in\RR_+}$ \ as \ $n\to\infty$, \ since in this case the process
 \ $(\cX_t)_{t\in\RR_+}$ \ given by \eqref{SDE_single} takes the form \ $\cX_t = m_\vare t$, \ $t\in\RR_+$.
\proofend
\end{Rem}

The next remark is devoted to a discussion on the role of the initial value \ $X_0$.

\begin{Rem}\label{Rem_initial_value}
Wei and Winnicki \cite[Theorem 2.1]{WW} do not suppose the finiteness of the second moment of the initial value \ $X_0$, \
 in their proof (which is based on infinitesimal generators) they only use that \ $n^{-1}X_0$ \ converges to \ $0$ \ as \ $n\to\infty$ \ almost surely,
 which holds without any further assumption on \ $X_0$.
In their proof, Wei and Winnicki \cite{WW} refer to several results of Ethier and Kurtz \cite{EthKur}, such as Theorem 1.3 in Chapter 9 and
 (implicitly) Theorem 6.5 in Chapter 1, Theorem 8.2 and Corollary 8.9 in Chapter 4,
 and one can see that the initial value comes into play in Theorem 8.2 in Chapter 4 in Ethier and Kurtz \cite{EthKur}.
Since we are not experts in the theory of infinitesimal generators, we can not give further insights into the role of
 the initial value in the proof of Theorem 2.1 in Wei and Winnicki \cite{WW}.
Next, we explain the role of the initial value \ $X_0$ \ in our proof of Theorem \ref{Thm_critical} and the second order moment assumption on it.
Note that \ $X_0$ \ appears in the definition of \ $\cM_t^{(n)}$ \ (see \eqref{helpM}),
 and one can also realize that one has to handle \ $X_0$ \  in proving \ $n^{-2} \sup_{t\in[0,T]} X_{\nt}\stoch 0$ \ as \ $n\to\infty$ \ for each \ $T\in\RR_{++}$
 \ (see \eqref{Cond11}).
Further, in the course of the proof of Theorem \ref{Thm_critical} we need some estimation for \ $\EE(X_k^2)$, \ $k\in\ZZ_+$,
 \ for which we need to assume that \ $\EE(X_0^2)<\infty$.
Such an estimation is presented in Lemma \ref{Lem_mom_asymp} which is based on Lemma \ref{Moments}, where
 explicit formulae are derived for the (conditional) first two moments of \ $X_k$, \ $k\in\NN$, \ and \ $M_k$, \ $k\in\NN$.
In fact, the proof of Lemma \ref{Moments} is presented mainly for instructional purposes in order to highlight the role
 of the initial value \ $X_0$.
\proofend
\end{Rem}

The paper is organized as follows.
Section \ref{Section_Proofs} is devoted to a detailed proof of Theorem \ref{Thm_critical} and to a comparison with the proof of Theorem 2.1 in Isp\'any \cite{Isp}.
We close the paper with three appendices: we recall formulae and estimates for
 first and second order moments of a critical Galton-Watson process with immigration (Appendix \ref{app:moments}),
 we recall a version of the continuous mapping theorem (Appendix \ref{CMT}) and a result about convergence of random step processes
 towards a diffusion process due to Isp\'any and Pap \cite{IspPap} (Appendix \ref{section_conv_step_processes}).

We decided to write this somewhat didactic note, since we wanted to understand clearly the role of the initial value \ $X_0$ \ in the proof of
 Theorem \ref{Thm_critical} and we wanted to present the usefulness of the limit theorem for random step processes (especially created from
 martingale differences) towards a diffusion process due to Isp\'any and Pap \cite{IspPap} directly in case of a critical Galton-Watson process
 with immigration instead of a specialization of a corresponding result or proof for more general branching processes.

\section{Proof of Theorem \ref{Thm_critical} and comparison with the proof of Theorem 2.1 in Isp\'any \cite{Isp}}
\label{Section_Proofs}

\noindent{\bf Proof of Theorem \ref{Thm_critical}.}
We divide the proof into several steps.
First, we prove weak convergence of a sequence of random step processes
 constructed from the martingale differences created from \ $(X_k)_{k\in\ZZ_+}$.
\ Namely, let us introduce the sequence
 \begin{equation*}
  M_k := X_k - \EE(X_k \mid \cF_{k-1}^{X})
        = X_k - X_{k-1} - m_\vare, \quad k \in \NN,
 \end{equation*}
 of martingale differences with respect to the filtration \ $(\cF_k^X)_{k\in\ZZ_+}$,
 where we used that \ $\EE(X_k \mid \cF_{k-1}^X) = X_{k-1}  + m_\vare$, \ $k\in\NN$, \ and recall
 that \ $m_\vare=\EE(\vare)$.
\ Consider the random step processes
 \begin{align}\label{helpM}
   \cM_t^{(n)} := \frac{1}{n} \biggl(X_0 + \sum_{k=1}^\nt M_k\biggr)
   = \frac{1}{n} X_\nt - \frac{\nt}{n} m_\vare , \qquad t \in \RR_+ , \qquad n \in \NN .
 \end{align}
We will show that
 \begin{equation}\label{conv_M}
  (\cM_t^{(n)})_{t\in\RR_+} \distr (\cM_t)_{t\in\RR_+} \qquad
  \text{as \ $n \to \infty$,}
 \end{equation}
 where the limit process \ $(\cM_t)_{t\in\RR_+}$ \ is the pathwise unique strong solution of the SDE
 \begin{equation}\label{SDE_M}
  \dd \cM_t
  = \sqrt{\sigma_\xi^2 (\cM_t + m_\vare t)^+ } \, \dd \cW_t , \qquad t\in\RR_+,
 \end{equation}
 with initial value \ $\cM_0 = 0$.
\ The proof of \eqref{conv_M} is based on a result due to Isp\'any and Pap \cite{IspPap} (see also Theorem \ref{Conv2DiffThm}),
 which is about convergence of random step processes towards a diffusion process.
Using weak convergence of \ $(\cM_t^{(n)})_{t\in\RR_+}$, \ an application of a version of the continuous mapping theorem
 (see Lemma \ref{lemma:kallenberg}) will yield weak convergence of \ $(n^{-1} X_\nt)_{t\in\RR_+}$ \ as \ $n\to\infty$.

{\sl Step 1 (starting steps for the proof of \eqref{conv_M}).}
\ In order to prove \eqref{conv_M}, we want to apply Theorem \ref{Conv2DiffThm} with \ $\cU := \cM$, \ $U_k^{(n)} := n^{-1} M_k$, $k\in\NN$,
 \ $U_0^{(n)} := n^{-1} X_0$, \ $\cF_k^{(n)} := \cF_k^X$, \ $k \in \ZZ_+$, \ where \ $n\in\NN$ \
 (yielding \ $\cU^{(n)} = \cM^{(n)}$, \ $n\in\NN$, \ as well),
 and with coefficient functions
 \ $\beta : \RR_+ \times \RR \to \RR$ \ and \ $\gamma : \RR_+ \times \RR \to \RR$ \ of the SDE \eqref{SDE_M} given by
 \[
   \beta(t, x) := 0 , \qquad
   \gamma(t, x) := \sqrt{\sigma_\xi^2 (x + m_\vare t)^+ } , \qquad t \in \RR_+ , \qquad x \in \RR .
 \]
First we check that the SDE \eqref{SDE_M} has a pathwise unique strong solution \ $(\cM_t^{(x)})_{t\in\RR_+}$ \ for all initial values \ $\cM_0^{(x)} = x \in \RR$.
\ Observe that if \ $(\cM_t^{(x)})_{t\in\RR_+}$ \ is a strong solution of the SDE \eqref{SDE_M} with initial value \ $\cM_0^{(x)} = x \in \RR$,
 \ then, by It\^o's formula, the process \ $\cP_t := \cM_t^{(x)} + m_\vare t$, \ $t \in \RR_+$, \ is a strong solution of the SDE
 \begin{equation}\label{SDE_P}
  \dd \cP_t
  = m_\vare \, \dd t
    + \sqrt{\sigma_\xi^2 \, \cP_t^+} \, \dd \cW_t, \qquad t\in\RR_+,
 \end{equation}
 with initial value \ $\cP_0 = x$.
\ Conversely, if \ $(\cP_t^{(p)})_{t\in\RR_+}$ \ is a strong solution of the SDE \eqref{SDE_P} with initial value \ $\cP_0^{(p)}  = p \in \RR$,
 \ then, by It\^o's formula, the process \ $\cM_t := \cP_t^{(p)} - m_\vare t$, \ $t \in \RR_+$, \ is a strong solution of the SDE \eqref{SDE_M} with initial value \ $\cM_0 = p$.
\ The SDE \eqref{SDE_P} is the same as \eqref{SDE_single}.
Consequently, as it was explained after Theorem \ref{Thm_critical}, the SDE \eqref{SDE_P} and hence the SDE \eqref{SDE_M} as well admit
 a pathwise unique strong solution with arbitrary initial value, and \ $(\cM_t + m_\vare t)_{t\in\RR_+} \distre (\cX_t)_{t\in\RR_+}$.

Note that \ $\EE( (U_k^{(n)})^2 )<\infty$ \ for all \ $n\in\NN$ \ and \ $k\in\ZZ_+$, \ since, by Lemma \ref{Moments},
 \ $\EE((U_k^{(n)})^2 ) = n^{-2}\EE(M_k^2)<\infty$, \ $n,k\in\NN$, \ and, by the assumption, \ $\EE((U_0^{(n)})^2 ) = n^{-2}\EE(X_0^2)<\infty$, \ $n\in\NN$.
\ Further, \ $U_0^{(n)} = n^{-1} X_0 \as 0$ \ as \ $n\to\infty$, \ especially \ $U_0^{(n)}\distr 0$ \ as \ $n\to\infty$.

For conditions (i), (ii) and (iii) of Theorem \ref{Conv2DiffThm}, we have to check that for each \ $T \in \RR_{++}$,
 \begin{gather} \label{Cond0}
  \sup_{t\in[0,T]}
   \biggl|\frac{1}{n}
          \sum_{k=1}^\nt
           \EE(M_k \mid \cF_{k-1}^X) - 0\biggr|
     \stoch 0\qquad \text{as \ $n\to\infty$,}\\
   \label{Cond1}
  \sup_{t\in[0,T]}
   \biggl|\frac{1}{n^2}
          \sum_{k=1}^\nt
           \EE(M_k^2 \mid \cF_{k-1}^X)
          - \int_0^t \sigma_\xi^2 (\cM_s^{(n)} + m_\vare s)^+  \, \dd s\biggr|
  \stoch 0\qquad \text{as \ $n\to\infty$,} \\
  \frac{1}{n^2}
  \sum_{k=1}^\nT
   \EE(M_k^2 \bbone_{\{|M_k|>n\theta\}} \mid \cF_{k-1}^X)
  \stoch 0  \qquad \text{as \ $n\to\infty$ \ for all \ $\theta \in \RR_{++}$.} \label{Cond2}
 \end{gather}
Condition \eqref{Cond0} trivially holds, since \ $\EE(M_k \mid \cF_{k-1}^X) =0$, \ $n\in\NN$, \ $k\in\NN$.

{\sl Step 2 (checking \eqref{Cond1}).}
For each \ $s \in \RR_+$ \ and \ $n \in \NN$, \ we have
 \[
   \cM_s^{(n)} + m_\vare s
   = \frac{1}{n} X_\ns + \frac{ns-\ns}{n} m_\vare ,
 \]
 thus \ $(\cM_s^{(n)} + m_\vare s)^+ = \cM_s^{(n)} + m_\vare s$, \ and
 \begin{align*}
  &\int_0^t (\cM_s^{(n)} + m_\vare s)^+ \, \dd s
   = \int_0^t \biggl(\frac{1}{n} X_\ns + \frac{ns-\ns}{n} m_\vare\biggr) \dd s \\
  &= \sum_{k=0}^{\nt-1} \int_{k/n}^{(k+1)/n} \biggl(\frac{1}{n} X_k + \frac{ns-k}{n} m_\vare\biggr) \dd s
     + \int_{\nt/n}^t \biggl(\frac{1}{n} X_\nt + \frac{ns-\nt}{n} m_\vare\biggr) \dd s \\
  &= \frac{1}{n^2} \sum_{k=0}^{\nt-1} X_k
     + \frac{m_\vare}{n}
       \sum_{k=0}^{\nt-1} \Bigl[\frac{1}{2} n s^2 - k s\Bigr]_{s=k/n}^{s=(k+1)/n} \\
  &\quad
     + \frac{1}{n} \biggl(t - \frac{\nt}{n}\biggr) X_\nt
     + \frac{m_\vare}{n} \Bigl[\frac{1}{2} n s^2 - \nt s\Bigr]_{s=\nt/n}^{s=t}\\
  &= \frac{1}{n^2} \sum_{k=0}^{\nt-1} X_k
     + \frac{nt-\nt}{n^2} X_\nt
     + \frac{m_\vare}{n} \left( \frac{n}{2}\cdot \frac{\nt^2 }{n^2} - \frac{\nt (\nt -1)}{2n} \right)\\
  &\quad + \frac{m_\vare}{n} \left( \frac{n}{2}\left( t^2 - \frac{\nt^2}{n^2}\right) - \nt \left(t-\frac{\nt}{n}\right) \right)\\
  & = \frac{1}{n^2} \sum_{k=0}^{\nt-1} X_k
     + \frac{nt-\nt}{n^2} X_\nt
     + \frac{\nt+(nt-\nt)^2}{2n^2} m_\vare
 \end{align*}
 for all \ $t \in \RR_+$ \ and \  $n \in \NN$.
\ By \eqref{condvar_X},
 \[
   \frac{1}{n^2}
   \sum_{k=1}^\nt
    \EE(M_k^2 \mid \cF_{k-1}^X)
   = \frac{\nt}{n^2} \sigma_\vare^2 + \frac{\sigma_\xi^2}{n^2} \sum_{k=1}^\nt  X_{k-1} , \qquad t \in \RR_+ , \qquad n \in \NN ,
 \]
 which yields that
 \begin{align*}
   &\frac{1}{n^2} \sum_{k=1}^\nt \EE(M_k^2 \mid \cF_{k-1}^X)
     - \int_0^t \sigma_\xi^2 (\cM_s^{(n)} + m_\vare s)^+ \, \dd s\\
   &\qquad = \frac{\nt}{n^2} \sigma_\vare^2
      - \sigma_\xi^2 \frac{nt-\nt}{n^2} X_{\nt}
      - \sigma_\xi^2 m_\vare \frac{\nt + (nt-\nt)^2}{2n^2},
      \qquad t\in\RR_+, \quad n\in\NN.
 \end{align*}
Since for each \ $T\in\RR_{++}$,
 \begin{align*}
  &\sup_{t\in[0,T]} \frac{\nt}{n^2} \leq \frac{T}{n}\to 0 \qquad \text{as \ $n\to\infty$,}\\
  &\sup_{t\in[0,T]}  \frac{\nt + (nt-\nt)^2}{2n^2} \leq \frac{T}{2n} + \frac{1}{2n^2}\to 0 \qquad \text{as \ $n\to\infty$,}
 \end{align*}
 in order to show \eqref{Cond1}, it suffices to prove that for each \ $T\in\RR_{++}$,
 \begin{equation}\label{Cond11}
  \frac{1}{n^2}
  \sup_{t \in [0,T]}
        \big((nt-\nt) X_\nt \big) \leq
  \frac{1}{n^2}
  \sup_{t \in [0,T]}
   X_\nt
  \stoch 0 \qquad \text{as \ $n \to \infty$.}
 \end{equation}
For each \ $k \in \NN$, \ we have \ $X_k = X_{k-1} + M_k + m_\vare$, \ thus
 \begin{equation*}
  X_k = X_0 + \sum_{j=1}^k M_j + k m_\vare ,
 \end{equation*}
 hence, for each \ $t \in \RR_+$ \ and \  $n \in \NN$, \ we get
 \[
   X_\nt = \vert X_\nt \vert \leq X_0 + \sum_{j=1}^\nt |M_j| + \nt m_\vare.
 \]
Consequently, in order to prove \eqref{Cond11}, it suffices to show
 \[
   \frac{1}{n^2}\sup_{t\in[0,T]} \sum_{j=1}^\nt |M_j|
    \leq \frac{1}{n^2} \sum_{j=1}^\nT |M_j| \stoch 0 \qquad \text{as \ $n \to \infty$.}
 \]
By Lemma \ref{Lem_mom_asymp},
 \[
   \EE\Biggl(\frac{1}{n^2} \sum_{j=1}^\nT |M_j|\Biggr)
   = \frac{1}{n^2} \sum_{j=1}^\nT \OO(j^{1/2})
   = \OO(n^{-1/2})\to 0 \qquad \text{as \ $n\to\infty$,}
 \]
 thus we obtain \ $n^{-2}\sum_{j=1}^{\lfloor nT\rfloor}\vert M_j\vert\stoch 0$ \ as \ $n\to\infty$ \ yielding
 \eqref{Cond11}, and hence \eqref{Cond1}, as desired.

{\sl Step 3 (checking \eqref{Cond2}).}
In order to prove \eqref{Cond2}, for each \ $k \in \NN$, \ consider the decomposition
  of \ $M_k$ \ into a random sum of independent centered random variables and another centered random variable, which are independent, namely,
 \[
   M_k = \sum_{j=1}^{X_{k-1}} \xi_{k,j} + \vare_k - X_{k-1} - m_\vare
   = N_k + (\vare_k - m_\vare)
 \]
 with
 \[
   N_k := \sum_{j=1}^{X_{k-1}} (\xi_{k,j} - 1) .
 \]
For each \ $n,k \in \NN$ \ and \ $\theta \in \RR_{++}$, \ we have
 \[
   M_k^2 \leq 2 (N_k^2 + (\vare_k - m_\vare)^2) , \qquad
   \bbone_{\{|M_k|>n\theta\}}
   \leq \bbone_{\{|N_k|>n\theta/2\}} + \bbone_{\{|\vare_k - m_\vare|>n\theta/2\}} ,
 \]
 yielding
 \begin{align*}
   M_k^2 \bbone_{\{|M_k|>n\theta\}}
     \leq 2N_k^2 \bbone_{\{|N_k|>n\theta/2\}}
          + 2N_k^2 \bbone_{\{|\vare_k - m_\vare|>n\theta/2\}}
          + 2 (\vare_k - m_\vare)^2,
 \end{align*}
 and hence \eqref{Cond2} will be proved once we show
 \begin{gather}
  \frac{1}{n^2}
   \sum_{k=1}^{\nT}
    \EE(N_k^2 \bbone_{\{|N_k|>n\theta\}} \mid \cF_{k-1}^X)
   \stoch 0 \qquad \text{as \ $n\to\infty$ \ for all \ $\theta \in \RR_{++}$,} \label{Cond21} \\
  \frac{1}{n^2}
   \sum_{k=1}^{\nT}
    \EE(N_k^2 \bbone_{\{|\vare_k-m_\vare|>n\theta\}}
        \mid \cF_{k-1}^X)
   \stoch 0 \qquad \text{as \ $n\to\infty$ \ for all \ $\theta \in \RR_{++}$,} \label{Cond22} \\
  \frac{1}{n^2}
   \sum_{k=1}^{\nT}
    \EE((\vare_k - m_\vare)^2 \mid \cF_{k-1}^X)
   \stoch 0 \qquad \text{as \ $n\to\infty$.} \label{Cond23}
 \end{gather}
In what follows let \ $\theta\in\RR_{++}$ \ be fixed.

{\sl Step 3/a (checking  \eqref{Cond21}).}
Using that the random variables \ $\{\xi_{k,j} : j \in \NN\}$
 \ are independent of the \ $\sigma$-algebra \ $\cF_{k-1}^X$ \ for all \ $k \in \NN$, \
 by the properties of conditional expectation with respect to a \ $\sigma$-algebra,
 we get for all \ $n,k\in\NN$,
 \[
   \EE(N_k^2 \bbone_{\{|N_k|>n\theta\}} \mid \cF_{k-1}^X)
   = F_{n,k}(X_{k-1}) ,
 \]
 where \ $F_{n,k} : \ZZ_+ \to \RR$ \ is given by
 \[
   F_{n,k}(z)
   :=\EE(S_k(z)^2 \bbone_{\{ \vert S_k(z)\vert >n\theta\}}) ,
   \qquad z \in \ZZ_+ ,
 \]
 with
 \begin{align}\label{help_S_k}
   S_k(z) :=  \sum_{j=1}^z (\xi_{k,j} - 1) , \qquad z \in \ZZ_+ .
 \end{align}
Consider the decomposition \ $F_{n,k}(z) = A_{n,k}(z) + B_{n,k}(z)$ \ with
 \begin{align*}
  &A_{n,k}(z) :=  \sum_{j=1}^z \EE\big( (\xi_{k,j} - 1)^2 ) \bbone_{\{\vert S_k(z)\vert >n\theta\}}  \big) ,  \\
  &B_{n,k}(z) :=  \sum\nolimits'_{j,j'} \EE\big( (\xi_{k,j} - 1) (\xi_{k,j'} - 1) ) \bbone_{\{\vert S_k(z)\vert>n\theta\}} \big),
 \end{align*}
 where the sum \ $\sum'_{j,j'}$ \ is taken for \ $j, j' \in \{1, \dots, z\}$ \ with \ $j \ne j'$.
\ Consider the inequalities
 \[
   \vert S_k(z)\vert = \vert \xi_{k,j} - 1 + \tS_k^j(z) \vert  \leq |\xi_{k,j} - 1| + \vert \tS_k^j(z)\vert, \qquad z\in\ZZ_+,
 \]
 for \ $j \in \{1, \ldots, z\}$, \ where
 \[
   \tS_k^j(z) := \sum\nolimits''_{j'} (\xi_{k,j'} - 1) ,\qquad z\in\ZZ_+,
 \]
 where the sum \ $\sum''_{j'}$ \ is taken for \ $j' \in \{1, \ldots, z\}$ \ with \ $j' \ne j$.
\ Using that
 \[
   \bbone_{\{\vert S_k(z)\vert >n\theta\}}
   \leq \bbone_{\{| \xi_{k,j} - 1|>n\theta/2\}} + \bbone_{ \{\vert \tS^j_k(z)\vert >n\theta/2\} } ,
   \qquad j\in\{1,\ldots,z\},
 \]
 we have \ $A_{n,k}(z) \leq A_{n,k}^{(1)}(z) + A_{n,k}^{(2)}(z)$, \ where
 \begin{align*}
  A_{n,k}^{(1)}(z) &:= \sum_{j=1}^z \EE( (\xi_{k,j} - 1)^2   \bbone_{\{|\xi_{k,j} - 1|>n\theta/2\}} ), \\
  A_{n,k}^{(2)}(z) &:= \sum_{j=1}^z \EE( (\xi_{k,j} - 1)^2 \bbone_{\{\vert \tS^j_k(z)\vert >n\theta/2\}} ) .
 \end{align*}
In order to prove \eqref{Cond21}, it is enough to show that
 \begin{gather}\nonumber
  \frac{1}{n^2} \sum_{k=1}^\nT A_{n,k}^{(1)}(X_{k-1})
  \stoch 0 , \qquad \qquad
  \frac{1}{n^2} \sum_{k=1}^\nT A_{n,k}^{(2)}(X_{k-1})
  \stoch 0 , \\ \label{Bk}
  \frac{1}{n^2} \sum_{k=1}^\nT B_{n,k}(X_{k-1}) \stoch 0
 \end{gather}
 as \ $n \to \infty$.
\ Using that \ $\xi_{k,j}$, $k,j\in\NN$, \ are identically distributed we have
 \[
   A_{n,k}^{(1)}(z) = z \EE( (\xi_{1,1} - 1)^2 \bbone_{\{|\xi_{1,1} - 1|>n\theta/2\}}) , \qquad n,k \in \NN ,\quad z\in\ZZ_+,
 \]
 thus, by Lemma \ref{Lem_mom_asymp}, we get
 \begin{align*}
  &\EE\biggl(\frac{1}{n^2} \sum_{k=1}^\nT A_{n,k}^{(1)}(X_{k-1})\biggr)
   = \frac{1}{n^2}
     \left(\sum_{k=1}^\nT \EE(X_{k-1})\right)
     \EE((\xi_{1,1} - 1)^2 \bbone_{\{|\xi_{1,1}-1|>n\theta/2\}}) \\
  &= \frac{1}{n^2}
     \left(\sum_{k=1}^\nT \OO(k)\right) \EE((\xi_{1,1} - 1)^2 \bbone_{\{|\xi_{1,1}-1|>n\theta/2\}})
   = \EE((\xi_{1,1} - 1)^2 \bbone_{\{|\xi_{1,1}-1|>n\theta/2\}}) \OO(1)
 \end{align*}
 for \ $n \in \NN$.
\ Consequently, since \ $\EE(\xi_{1,1}^2)<\infty$, \ the dominated convergence theorem implies
 \ $\EE\bigl(n^{-2} \sum_{k=1}^{\nT} A_{n,k}^{(1)}(X_{k-1})\bigr) \to 0$ \ as \ $n \to \infty$, \ which yields
 \ $n^{-2} \sum_{k=1}^{\nT} A_{n,k}^{(1)}(X_{k-1}) \stoch 0$ \ as \ $n \to \infty$, \ as desired.

Further, the independence of \ $\xi_{k,j} - 1$ \ and \ $\tS^j_k(z)$ \ implies
 \[
   A_{n,k}^{(2)}(z) = \sum_{j=1}^z \EE((\xi_{k,j} - 1)^2) \PP\bigl( \vert \tS^j_k(z)\vert  > n\theta/2\bigr) ,
    \qquad n,k\in\NN, \quad z\in\ZZ_+.
 \]
Here \ $\EE((\xi_{k,j} - 1)^2) = \sigma_\xi^2$, \ and, by \ $\EE(\tS^j_k(z))=0$ \ and Markov inequality, we have
 \[
   \PP\bigl( \vert \tS^j_k(z) \vert> n\theta/2\bigr)
     = \PP\bigl( \vert \tS^j_k(z) - \EE(\tS^j_k(z)) \vert > n\theta/2\bigr)
   \leq  \frac{4}{n^2\theta^2} \var(\tS^j_k(z))
   = \frac{4}{n^2\theta^2} \EE(\tS^j_k(z)^2),
  \]
 where, using that \ $\xi_{1,1}$ \ and \ $\xi_{1,2}$ \ are independent,
 \begin{align*}
   \EE(\tS^j_k(z)^2)
     &= \EE \Big(\Big( \sum\nolimits''_{j'} (\xi_{k,j'} - 1) \Big)^2 \Big)
      = \EE \Big( \sum\nolimits''_{j'} (\xi_{k,j'} - 1) \sum\nolimits''_{\ell'} (\xi_{k,\ell'} - 1)   \Big)\\
     & = \sum\nolimits''_{j'} \EE( (\xi_{k,j'} - 1)^2 ) + \sum\nolimits''_{j',\ell', j'\ne \ell'} \EE( (\xi_{k,j'} - 1)(\xi_{k,\ell'} - 1) )\\
     &= (z-1) \var(\xi_{1,1}) + (z-1)(z-2) \EE( (\xi_{1,1} - 1)(\xi_{1,2} - 1) ) \\
     &= (z-1) \sigma_\xi^2 + (z-1)(z-2) \EE( \xi_{1,1} - 1 ) \EE(\xi_{1,2} - 1 )
     = (z-1) \sigma_\xi^2 \leq z \sigma_\xi^2,
 \end{align*}
 where \ $\sum\nolimits''_{j'} $ \ and \ $\sum\nolimits''_{\ell'}$ \ is taken for \ $j'\in\{1,\ldots,z\}$ \ with \ $j'\ne j$, \
 and \ $\ell'\in\{1,\ldots,z\}$ \ with \ $\ell'\ne j$, \ respectively, and
 \ $\sum\nolimits''_{j',\ell', j'\ne \ell'}$ \ is taken for \ $j',\ell'\in\{1,\ldots,z\}$ \ with \ $j'\ne j$, \ $\ell'\ne \ell$, \ $j'\ne \ell'$.
\ Hence
 \[
   \PP\bigl( \vert \tS^j_k(z) \vert > n\theta/2\bigr)
     \leq \frac{4}{n^2\theta^2} z \sigma_\xi^2, \qquad z\in\ZZ_+,\quad j\in\{1,\ldots,z\},
 \]
 and consequently
 \begin{align*}
   A_{n,k}^{(2)}(z) \leq \frac{4}{n^2\theta^2} z^2 \sigma_\xi^4 , \qquad n,k\in\NN, \quad z\in\ZZ_+.
 \end{align*}
Consequently, by Lemma \ref{Lem_mom_asymp},
 \[
   \EE\biggl(\frac{1}{n^2} \sum_{k=1}^\nT A_{n,k}^{(2)}(X_{k-1})\biggr)
   \leq \frac{4\sigma_\xi^4}{n^4\theta^2}
        \sum_{k=1}^{\nT} \EE(X_{k-1}^2)
   = \frac{4 \sigma_\xi^4 }{n^4\theta^2}
     \sum_{k=1}^{\nT} \OO(k^2)
   = \OO(n^{-1})
 \]
 for \ $n \in \NN$, \ which implies \ $\EE\bigl(n^{-2} \sum_{k=1}^{\nT} A_{n,k}^{(2)}(X_{k-1})\bigr) \to 0$ \ as \ $n \to \infty$, \ and hence
 \ $n^{-2} \sum_{k=1}^{\nT} A_{n,k}^{(2)}(X_{k-1}) \stoch 0$ \ as \ $n \to \infty$, \ as desired.

Now we turn to check \eqref{Bk}.
By Cauchy-Schwarz inequality,
 \[
   |B_{n,k}(z)| \leq  \EE \left( \left\vert \sum\nolimits'_{j,j'} (\xi_{k,j} - 1) (\xi_{k,j'} - 1) \right\vert
                                       \bbone_{\{ \vert S_k(z)\vert >n\theta\}}  \right)
                 \leq  \sqrt{B_{n,k}^{(1)}(z) \, \EE(\bbone_{\{ \vert S_k(z)\vert >n\theta\}})},
                 \quad z\in\ZZ_+,
 \]
 where
 \[
   B_{n,k}^{(1)}(z)
   := \EE\biggl(\biggl(\sum\nolimits'_{j,j'} (\xi_{k,j} - 1) (\xi_{k,j'} - 1)\biggr)^2\biggr) .
 \]
Using the independence of \ $\xi_{k,j} - 1$ \ and \ $\xi_{k,j'} - 1$ \ for \ $j \ne j'$, \
 and that \ $\EE(\xi_{1,1} - 1)=0$, \ we get
 \[
   B_{n,k}^{(1)}(z)
   = 2z (z - 1) \sigma_\xi^4 \leq 2 \sigma_\xi^4 z^2 , \qquad z\in\ZZ_+,\quad n,k\in\NN.
 \]
Indeed, for \ $z\in\ZZ_+$,
 \begin{align*}
 B_{n,k}^{(1)}(z)
  & = \EE \biggl( \sum\nolimits'_{j,j'} (\xi_{k,j} - 1) (\xi_{k,j'} - 1)
      \sum\nolimits'_{\ell,\ell'} (\xi_{k,\ell} - 1) (\xi_{k,\ell'} - 1) \biggr) \\
  & = \biggl( \sum_{j=\ell,j'=\ell', j\ne j'} +  \sum_{j=\ell',j'=\ell, j\ne j'} \biggr)
       \EE( (\xi_{k,j} - 1)^2)  \EE( (\xi_{k,j'} - 1)^2)
    = 2z (z - 1) \sigma_\xi^4,
 \end{align*}
 and both sums \ $\sum_{j=\ell,j'=\ell', j\ne j'}$ \ and \ $\sum_{j=\ell',j'=\ell, j\ne j'}$ \ have \ $z^2-z = z(z-1)$ \ terms.
Further, by \ $\EE(S_k(z))=0$, \ $z\in\ZZ_+$, \ Markov inequality, and the independence of \ $\xi_{k,j}$, \ $j\in\ZZ_+$, \ we have
 \begin{align*}
   \EE(\bbone_{\{ \vert S_k(z) \vert >n\theta\}})
     = \PP(\vert S_k(z) - \EE(S_k(z)) \vert >n\theta )
           \leq \frac{\var(S_k(z))}{n^2\theta^2}
    = \frac{1}{n^2\theta^2} z \sigma_\xi^2 .
 \end{align*}
Hence
 \[
   |B_{n,k}(z)|
   \leq \sqrt{ 2 \sigma_\xi^4 z^2 n^{-2} \theta^{-2} z \sigma_\xi^2}
   = \frac{ \sqrt{2} \sigma_\xi^3 }{\theta n} z^{3/2} , \qquad z\in\ZZ_+,\quad n,k\in\NN.
 \]
Thus, in order to show \eqref{Bk}, it suffices to prove
 \[
   n^{-3} \sum_{k=1}^{\nT} X_{k-1}^{3/2}
   \stoch 0 \qquad \text{as \ $n \to \infty$.}
 \]
By Lyapunov inequality, \ $(\EE(X_{k-1}^{3/2}))^{2/3}\leq (\EE(X_{k-1}^2))^{1/2}$, \ and, using Lemma \ref{Lem_mom_asymp}, we get
 \[
   \EE(X_{k-1}^{3/2})
   \leq \big(\EE(X_{k-1}^2)\big)^{3/4}
   = \bigl(\OO(k^2)\bigr)^{3/4}
   = \OO(k^{3/2}) \qquad \text{for \ $k \in \NN$,}
 \]
 hence
 \[
   \EE\biggl(n^{-3} \sum_{k=1}^{\nT} X_{k-1}^{3/2}\biggr)
   = n^{-3} \sum_{k=1}^{\nT} \OO(k^{3/2})
   = \OO(n^{-1/2}) \qquad \text{for \ $n \in \NN$.}
 \]
Consequently, we obtain \ $\EE\bigl(n^{-3} \sum_{k=1}^{\nT} X_{k-1}^{3/2}\bigr)) \to 0$ \ as \ $n \to \infty$, \ yielding
 \ $n^{-3} \sum_{k=1}^{\nT} X_{k-1}^{3/2} \stoch 0$ \ as \ $n \to \infty$, \ yielding \eqref{Bk}, as desired.
Thus we finished the proof of \eqref{Cond21}.

{\sl Step 3/b (checking \eqref{Cond22}).}
Using that for all \ $k \in \NN$, \ the random variables
 \ $\bigl\{\xi_{k,j}, \vare_k : j \in \NN\bigr\}$
 \ are independent of the \ $\sigma$-algebra \ $\cF_{k-1}^X$, \
 by the properties of conditional expectation with respect to a \ $\sigma$-algebra, we get
 \[
   \EE(N_k^2 \bbone_{\{|\vare_k- m_\vare|>n\theta\}} \mid \cF_{k-1}^X) = G_k(X_{k-1}) ,
 \]
 where \ $G_k : \ZZ_+ \to \RR$ \ is given by
 \[
   G_k(z) := \EE(S_k(z)^2 \bbone_{\{|\vare_k- m_\vare|>n\theta\}}) , \qquad z \in \ZZ_+ ,
 \]
 and \ $S_k(z)$ \ is given in \eqref{help_S_k}.
Using again the independence of \ $\bigl\{\xi_{k,j}, \vare_k : j \in \NN\bigr\}$ \ and that \ $\EE(\xi_{1,1}-1)=0$, \ we have
 \[
   G_k(z)
   = \PP(|\vare_k- m_\vare | > n \theta) \sum_{j=1}^z \EE((\xi_{k,j} - 1)^2) ,\qquad z\in\ZZ_+,
 \]
 where, by Markov inequality,
 \ $\PP(|\vare_k- m_\vare | > n\theta) \leq n^{-2} \theta^{-2} \EE((\vare_k- m_\vare)^2)
    = n^{-2} \theta^{-2} \sigma_\vare^2$
 \ and \ $\EE((\xi_{k,j} - 1)^2) = \sigma_\xi^2$.
\ Hence
 \[
  G_k(z) \leq n^{-2}\theta^{-2} \sigma_\vare^2 \sigma_\xi^2 z, \qquad z\in\ZZ_+,
 \]
 and, in order to show \eqref{Cond22}, it suffices to prove
 \[
   n^{-4} \sum_{k=1}^{\nT} X_{k-1} \stoch 0  \qquad \text{as \ $n\to\infty$.}
 \]
In fact, by Lemma \ref{Lem_mom_asymp},
 \ $\EE\bigl(n^{-4} \sum_{k=1}^{\nT} X_{k-1}\bigr) = n^{-4} \sum_{k=1}^{\nT} \OO(k) = \OO(n^{-2})$ \ for \ $n \in \NN$, \ implying
 \ $\EE\bigl(n^{-4} \sum_{k=1}^{\nT} X_{k-1}\bigr) \to 0$ \ as \ $n \to \infty$, \ and hence \ $n^{-4} \sum_{k=1}^{\nT} X_{k-1} \stoch 0$ \ as \ $n \to \infty$,
 \ as desired.

{\sl Step 3/c (checking \eqref{Cond23}).}
By the independence of \ $\vare_k$ \ and \ $\cF_{k-1}^X$,
 \[
   \frac{1}{n^2}
   \sum_{k=1}^{\nT} \EE((\vare_k-m_\vare)^2 \mid \cF_{k-1}^X)
   = \frac{1}{n^2} \sum_{k=1}^{\nT} \EE((\vare_k-m_\vare)^2)
   = \frac{\nT}{n^2} \sigma_\vare^2
   \to 0 \quad \text{as \ $n\to\infty$},
 \]
 thus we obtain \eqref{Cond23}.

By Steps 3/a, 3/b and 3/c, we get \eqref{Cond2}, and, by Theorem \ref{Conv2DiffThm}, we conclude convergence \eqref{conv_M}.

{\sl Step 4 (proof of \eqref{convX}).}
In order to prove convergence \eqref{convX}, we want to apply Lemma \ref{lemma:kallenberg} using \eqref{conv_M}.
For each \ $n \in \NN$, \ by \eqref{helpM}, \ $(n^{-1}X_{\nt})_{t\in\RR_+} = \Psi^{(n)}(\cM^{(n)})$, \
 where the mapping \ $\Psi^{(n)} : \DD(\RR_+, \RR) \to \DD(\RR_+, \RR)$ \ is given by
 \[
   (\Psi^{(n)}(f))(t)
    := f\biggl(\frac{\nt}{n}\biggr) + \frac{\nt}{n} m_\vare
 \]
 for \ $f \in \DD(\RR_+, \RR)$ \ and \ $t \in \RR_+$.
\ Indeed, by \eqref{helpM}, for all \ $n\in\NN$ \ and \ $t\in\RR_+$,
 \begin{align*}
  (\Psi^{(n)}(\cM^{(n)}))(t)
    = \cM^{(n)}_{\nt/n} + \frac{\nt}{n}m_\vare
    = \frac{1}{n} X_{\lfloor n\cdot \frac{\nt}{n} \rfloor} - \frac{\lfloor n\cdot \frac{\nt}{n} \rfloor }{n}m_\vare + \frac{\nt}{n}m_\vare
    = \frac{1}{n} X_\nt.
 \end{align*}
Further, by \eqref{SDE_P}, \ $\cX \distre \Psi(\cM)$, \ where the mapping
 \ $\Psi : \DD(\RR_+, \RR) \to \DD(\RR_+, \RR)$ \ is given by
 \[
   (\Psi(f))(t) := f(t) + m_\vare t , \qquad
   f \in \DD(\RR_+, \RR) , \qquad t \in \RR_+ .
 \]

{\sl Step 4/a (checking measurability of \ $\Psi^{(n)}$, \ $n\in\NN$, \ and \ $\Psi$).}
We can check the measurability of the mappings \ $\Psi^{(n)}$, \ $n \in \NN$, \ and \ $\Psi$
 \ similarly as in Barczy et al. \cite[page 603]{BarIspPap0}.
Continuity of \ $\Psi$ \ follows from the characterization of convergence in
 \ $\DD(\RR_+, \RR)$, \ see, e.g., Ethier and Kurtz
 \cite[Proposition 3.5.3]{EthKur}, thus we obtain the measurability of \ $\Psi$ \ as well.
For each \ $n \in \NN$, \ in order to prove measurability of \ $\Psi^{(n)}$, \ first we localize it.
For each \ $n, N \in \NN$, \ consider the stopped mapping
 \ $\Psi^{(n,N)} : \DD(\RR_+, \RR) \to \DD(\RR_+, \RR)$ \ given by
 \ $(\Psi^{(n,N)}(f))(t) := (\Psi^{(n)}(f))(t \land N)$ \ for
 \ $f \in \DD(\RR_+, \RR)$, \ $t \in \RR_+$.
\ For each \ $f \in \DD(\RR_+, \RR)$, \ $T \in \RR_{++}$ \ and \ $N \in [T, \infty)$, \ we have
 \ $(\Psi^{(n,N)}(f))(t) = (\Psi^{(n)}(f))(t)$, \ $t \in [0,T]$, \ hence
 \ $\sup_{t \in [0,T]} \vert(\Psi^{(n,N)}(f))(t) - (\Psi^{(n)}(f))(t)\vert \to 0$ \ as
 \ $N \to \infty$, \ and then \ $\Psi^{(n,N)}(f) \to \Psi^{(n)}(f)$ \ in \ $\DD(\RR_+, \RR)$
 \ as \ $N \to \infty$, \ see, e.g., Jacod and Shiryaev
 \cite[VI.1.17]{JacShi}.
Consequently, it suffices to show measurability of \ $\Psi^{(n,N)}$ \ for all
 \ $n, N \in \NN$.
\ We can write \ $\Psi^{(n,N)} = \Psi^{(n,N,2)} \circ \Psi^{(n,N,1)}$, \ where the
 mappings \ $\Psi^{(n,N,1)} : \DD(\RR_+, \RR) \to \RR^{nN+1}$ \ and
 \ $\Psi^{(n,N,2)} : \RR^{nN+1} \to \DD(\RR_+, \RR)$ \ are defined by
 \begin{align*}
  \Psi^{(n,N,1)}(f)
  &:= \biggl(f(0), f\biggl(\frac{1}{n}\biggr), f\biggl(\frac{2}{n}\biggr), \dots, f(N)\biggr) , \\
  (\Psi^{(n,N,2)}(x_0, x_1, \dots, x_{nN}))(t)
  &:= x_{\lfloor n(t \land N)\rfloor} + \frac{\lfloor n(t \land N)\rfloor}{n} m_\vare
 \end{align*}
 for \ $f \in \DD(\RR_+, \RR)$, \ $t \in \RR_+$,
 \ $x_0, x_1, \dots, x_{nN} \in \RR$,
 \ $n, N \in \NN$.
\ Measurability of \ $\Psi^{(n,N,1)}$ \ follows from Ethier and Kurtz
 \cite[Proposition 3.7.1]{EthKur}.
Next we show continuity of \ $\Psi^{(n,N,2)}$. \
By Jacod and Shiryaev \cite[VI.1.17]{JacShi}, it is enough to check that
 \ $\sup_{t \in [0,T]} | (\Psi^{(n,N,2)}(\bx^{(k)}))(t) - (\Psi^{(n,N,2)}(\bx))(t)| \to 0$
 \ as \ $k \to \infty$ \ for all \ $T \in \RR_{++}$ \ whenever \ $\bx^{(k)} = (x_0^{(k)}, x_1^{(k)}, \dots, x_{nN}^{(k)}) \to \bx = (x_0, x_1, \dots, x_{nN})$
 \ as \ $k\to\infty$ \ in \ $\RR^{nN+1}$.
 \ This convergence follows from the estimate
 \begin{align*}
  \sup_{t \in [0,T]} |(\Psi^{(n,N,2)}(\bx^{(k)}))(t) - (\Psi^{(n,N,2)} (\bx))(t)|
  &= \sup_{t \in [0,T]} |x^{(k)}_{\lfloor n(t \land N)\rfloor} - x_{\lfloor n(t \land N)\rfloor}| \\
  &= \max_{j\in\{0,1,\ldots,nN\}} |x^{(k)}_j - x_j|
   \leq \|\bx^{(k)} - \bx\| ,
 \end{align*}
 where \ $\Vert \cdot\Vert$ \ denotes Euclidean norm.
We obtain measurability of both \ $\Psi^{(n,N,1)}$ \ and \ $\Psi^{(n,N,2)}$,
 \ hence we conclude measurability of \ $\Psi^{(n,N)}$.

{\sl Step 4/b (checking condition of Lemma \ref{lemma:kallenberg}).}
The aim of the following discussion is to show that the set
  \ $C := \CC(\RR_+, \RR)$ \ satisfies \ $C \in \cB(\DD(\RR_+, \RR))$, \ $\PP(\cM \in C) = 1$, \ and
 \ $\Psi^{(n)}(f^{(n)}) \to \Psi(f)$ \ in \ $\DD(\RR_+, \RR)$ \ as \ $n \to \infty$ \ if \ $f^{(n)} \to f$ \ in \ $\DD(\RR_+, \RR)$ \ as \ $n \to \infty$
 \ with \ $f \in C$, \ $f^{(n)}\in \DD(\RR_+, \RR)$, \ $n\in\NN$.

First note that \ $\CC(\RR_+, \RR) \in \cB (\DD(\RR_+, \RR))$, \ see,
 e.g., Ethier and Kurtz \cite[Problem 3.11.25]{EthKur}.
In fact, the subset \ $\CC(\RR_+, \RR) \subset \DD(\RR_+, \RR)$ \ is closed, since its complement \ $\DD(\RR_+, \RR) \setminus \CC(\RR_+, \RR)$ \ is
 open.
Indeed, each function
 \ $f \in \DD(\RR_+, \RR) \setminus \CC(\RR_+, \RR)$ \ is
 discontinuous at some point \ $t_f \in \RR_+$, \ and, by the definition of the
 metric of \ $\DD(\RR_+, \RR)$, \ there exists \ $r_f \in\RR_{++}$ \ such that
 all \ $g \in \DD(\RR_+, \RR)$ \ is discontinuous at the point
 \ $t_f \in \RR_+$ \ whenever the distance of \ $g$ \ and \ $f$ \ is less than
 \ $r_f$.
\ Consequently, the set
 \ $\DD(\RR_+, \RR) \setminus \CC(\RR_+, \RR)$ \ is the union of
 open balls with center
 \ $f \in \DD(\RR_+, \RR) \setminus \CC(\RR_+, \RR)$ \ and radius
 \ $r_f$.

By the definition of a strong solution (see, e.g., Jacod and Shiryaev
 \cite[Definition 2.24, Chapter III]{JacShi}), \ $\cM$ \ has
 continuous sample paths almost surely, so we have \ $\PP(\cM \in C) = 1$.

Next, we fix a function \ $f \in C$ \ and a sequence \ $(f^{(n)})_{n\in\NN}$ \ in
 \ $\DD(\RR_+, \RR)$ \ with \ $f^{(n)} \to f$ \ in \ $\DD(\RR_+, \RR)$ \ as \ $n \to \infty$.
\ Then the continuity of \ $f$ \ implies \ $f^{(n)} \lu f$ \ as \ $n \to \infty$, \ see, e.g., Jacod and Shiryaev
 \cite[VI.1.17]{JacShi} (for the notation \ $\lu$, \ see the beginning of Appendix \ref{CMT}).
By the definition of \ $\Psi$, \ we have
 \ $\Psi(f) \in \CC(\RR_+, \RR)$.
\ Further, for each \ $n \in \NN$, \ we can write
 \[
   (\Psi^{(n)}(f^{(n)}))(t) = f^{(n)}\biggl(\frac{\nt}{n}\biggr) + \frac{\nt}{n} m_\vare , \qquad t \in \RR_+ ,
 \]
 hence we have for all \ $T\in\RR_{++}$ \ and \ $t\in[0,T]$,
 \begin{align*}
  |(\Psi^{(n)}(f^{(n)}))(t) - (\Psi(f))(t)|
  &\leq \biggl|f^{(n)}\biggl(\frac{\nt}{n}\biggr) - f(t)\biggr| + \frac{1}{n} m_\vare \\
  &\leq \biggl|f^{(n)}\biggl(\frac{\nt}{n}\biggr) - f\biggl(\frac{\nt}{n}\biggr)\biggr| + \biggl|f\biggl(\frac{\nt}{n}\biggr) - f(t)\biggr| + \frac{1}{n} m_\vare \\
  & \leq  \sup_{t\in[0,T]} |f^{(n)}(t) - f(t)| + \omega_T(f,n^{-1}) + \frac{1}{n} m_\vare ,
 \end{align*}
 where \ $\omega_T(f, \cdot)$ \ is the modulus of continuity of \ $f$ \ on
 \ $[0, T]$.
\ We have \ $\omega_T(f, n^{-1}) \to 0$ \ as \ $n \to \infty$ \ since \ $f$ \ is
 continuous (see, e.g., Jacod and Shiryaev \cite[VI.1.6]{JacShi}), and \ $\sup_{t\in[0,T]} |f^{(n)}(t) - f(t)| \to 0$ \ as \ $n \to \infty$, \ since \ $f^{(n)} \lu f$ \ as \ $n \to \infty$.
 \ Thus we conclude \ $\Psi^{(n)}(f^{(n)}) \lu \Psi(f)$ \ as \ $n \to \infty$, \ and hence,
  since \ $\Psi(f)\in \CC(\RR_+,\RR)$, \ we have \ $\Psi^{(n)}(f^{(n)}) \to \Psi(f)$ \ in \ $\DD(\RR_+, \RR)$ \ as \ $n \to \infty$,
 \ see, e.g., Jacod and Shiryaev \cite[VI.1.17]{JacShi}.

{\sl Step 4/c (application of Lemma \ref{lemma:kallenberg}).}
Using Steps 4/a and 4/b, we can apply Lemma \ref{lemma:kallenberg} and we obtain
 $(n^{-1} X_\nt)_{t\in\RR_+} = \Psi^{(n)}(\cM^{(n)}) \distr \Psi(\cM)$ as
  $n \to \infty$, where $( (\Psi(\cM))(t))_{t\in\RR_+}  = (\cM_t + m_\vare t)_{t\in\RR_+} \distre (\cX_t)_{t\in\RR_+}$
  (by It\^o's formula, see \eqref{SDE_P}).
\proofend

In the next remark we compare our proof of Theorem \ref{Thm_critical} with the proof of Theorem 2.1 in Isp\'any \cite{Isp}
 by pointing out the parts where we made some simplifications in the arguments.

\begin{Rem}\label{Rem_Isp_comparison}
Theorem \ref{Thm_critical} is a special case of Theorem 2.1 in Isp\'any \cite{Isp}
 by considering a single critical Galton--Watson process with immigration and by choosing \ $\alpha=0$ \ in Definition 1.1 in Isp\'any \cite{Isp}.
In our proof of Theorem \ref{Thm_critical} we follow the same procedure as in the one of Theorem 2.1 in Isp\'any \cite{Isp},
 namely, we also use Theorem \ref{Conv2DiffThm}, and note that equations (15) and (16) in Isp\'any \cite{Isp}
 correspond to our equations \eqref{Cond1} and \eqref{Cond2}, respectively.
In the course of the proof of \eqref{Cond1} we give an explicit expression for \ $\int_0^t (\cM^{(n)}_s + m_\vare s)^+ \,\dd s$ \
 (with the notations of \cite{Isp}, for \ $\int_0^t \widetilde \cN_+^{(n)}(s)\,\dd s$) \ and also for
 \ $n^{-2}\sum_{k=1}^{\nt} \EE(M_k^2 \mid \cF_{k-1}^X) - \int_0^t \sigma_\xi^2(\cM_s^{(n)}+m_\vare s)^+\,\dd s$,
 \ which explicit forms are not available in \cite{Isp}, and in this way we think that the proof of \eqref{Cond1} becomes more understandable.
Further, concerning the proof of \eqref{Cond1}, at some point one needs to check that \ $n^{-2}\sup_{t\in[0,T]} X_{\nt}\stoch 0$ \ as \ $n\to\infty$ \
 for each \ $T\in\RR_{++}$, \ and to do so we do not need to use Lyapunov and Cauchy--Schwarz inequalities compared to the proof of the corresponding formula
 on page 29 in Isp\'any \cite{Isp} (due to the facts that in our special case, with the notations of \cite{Isp},
 \ $\sigma_n^2=\sigma_\vare^2$, \ $n\in\NN$, \ and \ $m_n=1$, $n\in\NN$).
Finally, we mention that in Isp\'any \cite{Isp} it is only stated that (using the notations and numberings in \cite{Isp})
 the weak convergence in (14) yields \ $(n\sigma_n^2)^{-1}\cX^{(n)}\distr \cX$ \ as \ $n\to\infty$, \ but not detailed at all (see \cite[the top of page 32]{Isp}).
However, in Step 4 of our proof of Theorem \ref{Thm_critical} we give a detailed exposition of the above mentioned step in our special case
 using a version of the continuous mapping theorem due to Kallenberg (see Appendix \ref{CMT}).
\proofend
\end{Rem}

\vspace*{5mm}

\appendix

\vspace*{5mm}

\noindent{\bf\Large Appendices}

\section{Moment estimates for critical Galton-Watson processes with immigration}
\label{app:moments}

In the proofs we use some facts about the first and second order moments of the sequences \ $(X_k)_{k\in\ZZ_+}$ \
 and \ $(M_k)_{k\in\ZZ_+}$ \ in the critical case (i.e., when \ $m_\xi=1$).

\begin{Lem}\label{Moments}
Let \ $(X_k)_{k\in\ZZ_+}$ \ be a critical Galton-Watson process with immigration
 such that \ $\EE(X_0^2)<\infty$, \ $\EE(\xi^2) < \infty$ \ and \ $\EE(\vare^2) < \infty$.
\ Then for all \ $k \in \NN$ \ we have
 \begin{align}
  &\EE(X_k \mid \cF_{k-1}^X) = X_{k-1} + m_\vare, \label{condmean_X} \\
  &\EE(X_k)
   = \EE(X_{k-1}) + m_\vare
   = \EE(X_0) + m_\vare k, \label{mean_X}
 \end{align}
 \begin{align}
  &\var(X_k \mid \cF_{k-1}^X)
	  	= \EE(M_k^2 \mid \cF_{k-1}^X)
        = \var(M_k \mid \cF_{k-1}^X)
        = \sigma_\xi^2 X_{k-1} + \sigma_\vare^2 , \label{condvar_X}\\
  &\var(X_k) = \var(X_{k-1}) + \sigma_\xi^2 \EE(X_{k-1}) + \sigma_\vare^2 \label{var_X} \\
    &\phantom{\var(X_k)} = m_\vare \sigma_\xi^2 \frac{(k-1)k}{2} + (\sigma_\xi^2 \EE(X_0) + \sigma_\vare^2)k + \var(X_0), \nonumber \\
  &\EE(M_k) = 0, \label{mean_M}\\
  &\EE(M_k^2) = \sigma_\xi^2\EE(X_{k-1}) + \sigma_\vare^2
              = \sigma_\xi^2 m_\vare (k-1) + \sigma_\xi^2 \EE(X_0) +\sigma_\vare^2, \label{second_mean_M}
 \end{align}
 where we recall \ $m_\vare=\EE(\vare)$, \ $\sigma_\xi^2=\var(\xi)$ \ and \ $\sigma_\vare^2=\var(\vare)$.
\end{Lem}

We note that a version of Lemma \ref{Moments} for critical multi-type Galton-Watson processes with immigration starting from zero can be found
 in Isp\'any and Pap \cite[Lemma A.2]{IspPap2}, and for single type Galton-Watson processes with immigration starting from zero, see also Isp\'any \cite[Lemma 4.1]{Isp}.
In case of \ $X_0=0$, \ Lemma \ref{Moments} is a special case of the above mentioned results due to Isp\'any and Pap \cite{IspPap2} and Isp\'any \cite{Isp},
 respectively.
For completeness, we present a proof.

\noindent
\textbf{Proof of Lemma \ref{Moments}.}
We already checked \eqref{condmean_X} and \eqref{mean_X}, see \eqref{EXk}.
For all \ $k\in\NN$, \ we have
 \begin{align*}
   M_k = X_k - X_{k-1} - m_\vare
      = \sum_{j=1}^{X_{k-1}} \xi_{k,j} + \vare_k
        - \sum_{j=1}^{X_{k-1}} 1 - m_\vare
      = \sum_{j=1}^{X_{k-1}} (\xi_{k,j} - 1)
        + (\vare_k - m_\vare),
 \end{align*}
 and using the independence of \ $\xi_{k,j}$, $j\in\NN$, \ and \ $\vare_k$, \
 and that they are independent of \ $\cF_{k-1}^X$, \ it yields that
 \begin{align*}
   \EE(M_k^2 \mid \cF_{k-1}^X)
     & \! = \! \EE\left( \sum_{j=1}^{X_{k-1}} \sum_{\ell=1}^{X_{k-1}}  (\xi_{k,j} - 1) (\xi_{k,\ell} - 1)
                  + 2\sum_{j=1}^{X_{k-1}}(\xi_{k,j} - 1) (\vare_k -m_\vare)
                  + (\vare_k -m_\vare)^2 \,\Big \vert\, \cF_{k-1}^X
           \right) \\
     & = \sum_{j=1}^{X_{k-1}} \EE((\xi -1)^2) + \EE((\vare_k -m_\vare)^2)
       = \sigma_\xi^2 X_{k-1} + \sigma_\vare^2,
 \end{align*}
 implying \eqref{condvar_X}.

Now we turn to prove \eqref{var_X}.
By the law of total variance, \eqref{condmean_X} and \eqref{condvar_X}, we have
 \begin{align*}
   \var(X_k)
     & = \EE(\var(X_k \mid \cF_{k-1}^X)) + \var(\EE(X_k \mid \cF_{k-1}^X))
       = \EE( \sigma_\xi^2 X_{k-1} + \sigma_\vare^2 ) + \var(X_{k-1} + m_\vare) \\
     & = \sigma_\xi^2 \EE(X_{k-1}) + \sigma_\vare^2 + \var(X_{k-1}),
     \qquad k\in\NN.
 \end{align*}
Hence, using also \eqref{mean_X}, for all \ $k\in\NN$, \ we have
 \begin{align*}
   \begin{bmatrix}
     \EE(X_k) \\
     \var(X_k) \\
   \end{bmatrix}
  & = \begin{bmatrix}
       1 & 0 \\
       \sigma_\xi^2 & 1 \\
     \end{bmatrix}
     \begin{bmatrix}
     \EE(X_{k-1}) \\
     \var(X_{k-1}) \\
    \end{bmatrix}
    + \begin{bmatrix}
       m_\vare \\
       \sigma_\vare^2 \\
      \end{bmatrix}\\
  & = \begin{bmatrix}
       1 & 0 \\
       \sigma_\xi^2 & 1 \\
     \end{bmatrix}^k
     \begin{bmatrix}
       \EE(X_0) \\
       \var(X_0) \\
      \end{bmatrix}
     + \sum_{j=0}^{k-1}
        \begin{bmatrix}
       1 & 0 \\
       \sigma_\xi^2 & 1 \\
     \end{bmatrix}^j
     \begin{bmatrix}
       m_\vare \\
       \sigma_\vare^2 \\
      \end{bmatrix}\\
  & = \begin{bmatrix}
       1 & 0 \\
       k\sigma_\xi^2 & 1 \\
     \end{bmatrix}
      \begin{bmatrix}
       \EE(X_0) \\
       \var(X_0) \\
      \end{bmatrix}
    + \sum_{j=0}^{k-1}
       \begin{bmatrix}
       1 & 0 \\
       j\sigma_\xi^2 & 1 \\
      \end{bmatrix}
     \begin{bmatrix}
       m_\vare \\
       \sigma_\vare^2 \\
      \end{bmatrix}.
 \end{align*}
Consequently,
 \begin{align*}
 \var(X_k)
  = \sigma_\xi^2 \EE(X_0) k + \var(X_0) + \sum_{j=0}^{k-1} (m_\vare \sigma_\xi^2 j + \sigma_\vare^2),
  \qquad k\in\NN,
 \end{align*}
 yielding \eqref{var_X}.
Finally, \eqref{mean_M} and \eqref{second_mean_M} follow by \eqref{condmean_X}, \eqref{mean_X} and \eqref{condvar_X}.
\proofend

\begin{Lem}\label{Lem_mom_asymp}
Under the conditions of Lemma \ref{Moments}, we have
 \begin{align*}
   \EE(X_k) = \OO(k) , \qquad \EE(X_k^2) = \OO(k^2), \qquad \EE(|M_k|) = \OO(k^{1/2}), \qquad \EE(M_k^2) = \OO(k), \qquad k\in\NN.
 \end{align*}
\end{Lem}

\noindent
\textbf{Proof.}
It follows by \eqref{mean_X}, \eqref{var_X} and \eqref{second_mean_M} together with \ $\EE(\vert M_k\vert)\leq \sqrt{\EE(M_k^2)}$, \ $k\in\NN$.
\proofend

\section{A version of the continuous mapping theorem}
\label{CMT}

A function \ $f : \RR_+ \to \RR$ \ is called \emph{c\`adl\`ag} if it is right
 continuous with left limits.
\ Let \ $\DD(\RR_+, \RR)$ \ and \ $\CC(\RR_+, \RR)$ \ denote the space of
 all real-valued c\`adl\`ag and continuous functions on \ $\RR_+$, \ respectively.
Let \ $\cB(\DD(\RR_+, \RR))$ \ denote the Borel $\sigma$-algebra on \ $\DD(\RR_+, \RR)$ \ for the metric defined in Jacod and Shiryaev
 \cite[Chapter VI, (1.26)]{JacShi} (with this metric \ $\DD(\RR_+, \RR)$ \ is a complete and separable metric space and
 the topology induced by this metric is the so-called Skorokhod topology).
For a function \ $f \in \DD(\RR_+, \RR)$ \ and for a sequence \ $(f_n)_{n\in\NN}$
 \ in \ $\DD(\RR_+, \RR)$, \ we write \ $f_n \lu f$ \ if \ $(f_n)_{n\in\NN}$
 \ converges to \ $f$ \ locally uniformly, i.e., if
 \ $\sup_{t\in[0,T]} |f_n(t) - f(t)| \to 0$ \ as \ $n \to \infty$ \ for all
 \ $T > 0$.
For real-valued stochastic processes \ $(\cY_t)_{t \in \RR_+}$ \ and
 \ $(\cY^{(n)}_t)_{t \in \RR_+}$, \ $n \in \NN$, \ with c\`adl\`ag paths we write
 \ $\cY^{(n)} \distr \cY$ \ if the distribution of \ $\cY^{(n)}$ \ on the
 space \ $(\DD(\RR_+, \RR), \cB(\DD(\RR_+, \RR)))$ \ converges weakly to the
 distribution of \ $\cY$ \ on the space \ $(\DD(\RR_+, \RR), \cB(\DD(\RR_+, \RR)))$ \ as \ $n \to \infty$.
\ Equality in distribution is denoted by \ $\distre$.
\ If \ $\xi$ \ and \ $\xi_n$, \ $n \in \NN$, \ are random elements with values in a metric space \ $(E, d)$,
 \ then we denote by \ $\xi_n \distr \xi$ \ the weak convergence of the
 distribution of \ $\xi_n$ \ on the space \ $(E, \cB(E))$ \ towards the
 distribution of \ $\xi$ \ on the space \ $(E, \cB(E))$ \ as \ $n \to \infty$,
 \ where \ $\cB(E)$ \ denotes the Borel \ $\sigma$-algebra on \ $E$ \ induced by
 the given metric \ $d$.

The following version of the continuous mapping theorem can be found for
 example in Theorem 3.27 of Kallenberg \cite{Kal}.

\begin{Lem}\label{lemma:kallenberg}
Let \ $(S, d_S)$ \ and \ $(T, d_T)$ \ be metric spaces and
 \ $(\xi_n)_{n \in \NN}$, \ $\xi$ \ be random elements with values in \ $S$
 \ such that \ $\xi_n \distr \xi$ \ as \ $n \to \infty$.
\ Let \ $f : S \to T$ \ and \ $f_n : S \to T$, \ $n \in \NN$, \ be measurable
 mappings and \ $C \in \cB(S)$ \ such that \ $\PP(\xi \in C) = 1$ \ and
 \ $\lim_{n \to \infty} d_T(f_n(s_n), f(s)) = 0$ \ if
 \ $\lim_{n \to \infty} d_S(s_n,s) = 0$ \ and \ $s \in C$, \ $s_n\in S$, \ $n\in\NN$.
\ Then \ $f_n(\xi_n) \distr f(\xi)$ \ as \ $n \to \infty$.
\end{Lem}

\section{Convergence of random step processes}
\label{section_conv_step_processes}

We recall a result about convergence of one-dimensional random step processes towards a diffusion process,
 see Isp\'any and Pap \cite{IspPap}.

\begin{Thm}\label{Conv2DiffThm}
Let \ $\beta : \RR_+ \times \RR \to \RR$ \ and \ $\gamma : \RR_+ \times \RR \to \RR$ \ be continuous functions.
Assume that uniqueness in the sense of probability law holds for the SDE
 \begin{equation}\label{SDE}
  \dd \, \cU_t
  = \beta (t, \cU_t) \, \dd t + \gamma (t, \cU_t) \, \dd \cW_t ,
  \qquad t \in \RR_+,
 \end{equation}
 with initial value \ $\cU_0 = u_0$ \ for all \ $u_0 \in \RR$, \ where
 \ $(\cW_t)_{t\in\RR_+}$ \ is an one-dimensional standard Wiener process.
Let \ $(\cU_t)_{t\in\RR_+}$ \ be a solution of \eqref{SDE} with initial value
 \ $\cU_0 = 0$.

For each \ $n \in \NN$, \ let \ $(U^{(n)}_k)_{k\in\ZZ_+}$ \ be a sequence of
 real-valued random variables adapted to a filtration \ $(\cF^{(n)}_k)_{k\in\ZZ_+}$ \
 (i.e., \ $U^{(n)}_k$ \ is \ $\cF^{(n)}_k$-measurable) such that \ $\EE( (U^{(n)}_k)^2 )<\infty$ \ for each \ $n,k\in \NN$.
\ Let
 \[
   \cU^{(n)}_t := \sum_{k=0}^{\nt}  U^{(n)}_k \, ,
   \qquad t \in \RR_+, \quad n \in \NN .
 \]
Suppose that \ $\cU^{(n)}_0 = U^{(n)}_0 \distr 0$ \ as \ $n\to\infty$ \ and that for each \ $T \in \RR_{++}$,
 \begin{enumerate}
  \item[\textup{(i)}]
   $\sup\limits_{t\in[0,T]}
     \biggl\vert \sum\limits_{k=1}^{\nt}
              \EE\bigl(U^{(n)}_k \mid \cF^{(n)}_{k-1}\bigr)
             - \int_0^t \beta(s,\cU^{(n)}_s) \dd s\biggr\vert
   \stoch 0$  \ as \ $n\to\infty$, \\
  \item[\textup{(ii)}]
   $\sup\limits_{t\in[0,T]}
     \biggl\vert \sum\limits_{k=1}^{\nt}
              \var\bigl(U^{(n)}_k \mid \cF^{(n)}_{k-1}\bigr)
             - \int_0^t
                (\gamma(s,\cU^{(n)}_s) )^2
                \dd s\biggr\vert
   \stoch 0$ \ as \ $n\to\infty$, \\
  \item[\textup{(iii)}]
   $\sum\limits_{k=1}^{\lfloor nT \rfloor}
     \EE\bigl( (U^{(n)}_k)^2 \bbone_{\{\vert U^{(n)}_k\vert > \theta\}}
             \bmid \cF^{(n)}_{k-1}\bigr)
   \stoch 0$ \ as \ $n\to\infty$ \ for all \ $\theta \in \RR_{++}$.
 \end{enumerate}
Then \ $\cU^{(n)} \distr \cU$ \ as \ $n \to \infty$.
\end{Thm}

\section*{Acknowledgements}
This note was finished after the sudden death of the third author Gyula Pap in October 2019.
He liked the applications of Theorem \ref{Conv2DiffThm} in the theory of branching processes very much.
We would like to thank the referee for the suggestions that helped us to improve the paper.

\end{document}